\documentclass[12pt,a4paper]{article}
\usepackage{epsfig,latexsym,amsfonts,amssymb,amsmath,amscd,
graphics,theorem,epic}
\theoremstyle{plain}
\newtheorem{lemma}{Lemma}

\newtheorem{theorem}{Theorem}

{\theorembodyfont{\rmfamily} 
\font\ncsc=cmcsc10  \font\ntt=cmtt12
\begin{document}
\baselineskip=15pt
\newcommand{\pperp}{\hbox{$\perp\hskip-6pt\perp$}}
\newcommand{\ssim}{\hbox{$\hskip-2pt\sim$}}
\newcommand{\N}{{\mathbb N}}\newcommand{\Delp}{{\Pi}}
\newcommand{\A}{{\mathbb A}}
\newcommand{\Z}{{\mathbb Z}}
\newcommand{\R}{{\mathbb R}}
\newcommand{\C}{{\mathbb C}}
\newcommand{\Q}{{\mathbb Q}}\newcommand{\T}{{\mathbb T}}\newcommand{\K}{{\mathbb K}}
\newcommand{\PP}{{\mathbb P}}
\newcommand{\bbS}{{\mathbb S}}
\newcommand{\st}{{*}}
\newcommand{\mnote}{\marginpar}\newcommand{\ev}{{\operatorname{ev}}}
\newcommand{\Id}{{\operatorname{Id}}}\newcommand{\irr}{{\operatorname{irr}}}\newcommand{\nod}{{\operatorname{nod}}}
\newcommand{\oeps}{{\overline\eps}}\newcommand{\Area}{{\operatorname{Area}}}\newcommand{\End}{{\operatorname{End}}}
\newcommand{\oDel}{{\widetilde\Del}}
\newcommand{\real}{{\operatorname{Re}}}
\newcommand{\conv}{{\operatorname{conv}}}
\newcommand{\Span}{{\operatorname{Span}}}
\newcommand{\Ker}{{\operatorname{Ker}}}
\newcommand{\Fix}{{\operatorname{Fix}}}
\newcommand{\sign}{{\operatorname{sign}}}
\newcommand{\Log}{{\operatorname{Log}}}
\newcommand{\oi}{{\overline i}}
\newcommand{\oj}{{\overline j}}
\newcommand{\ob}{{\overline b}}
\newcommand{\os}{{\overline s}}
\newcommand{\oa}{{\overline a}}
\newcommand{\oy}{{\overline y}}
\newcommand{\ow}{{\overline w}}
\newcommand{\ou}{{\overline u}}
\newcommand{\ot}{{\overline t}}
\newcommand{\oz}{{\overline z}}
\newcommand{\newi}{i}
\newcommand{\newj}{j}
\newcommand{\newm}{m}
\newcommand{\newl}{{\ell}}
\newcommand{\bw}{{\boldsymbol w}}\newcommand{\bi}{{\boldsymbol i}}
\newcommand{\bx}{{\boldsymbol p}}\newcommand{\bp}{{\boldsymbol p}}
\newcommand{\bpp}{{\boldsymbol P}}
\newcommand{\by}{{\boldsymbol q}}
\newcommand{\bz}{{\boldsymbol z}}
\newcommand{\eps}{{\varepsilon}}
\newcommand{\proofend}{\hfill$\Box$\bigskip}
\newcommand{\Int}{{\operatorname{Int}}}
\newcommand{\pr}{{\operatorname{pr}}}
\newcommand{\grad}{{\operatorname{grad}}}
\newcommand{\rk}{{\operatorname{rk}}}
\newcommand{\im}{{\operatorname{Im}}}
\newcommand{\sk}{{\operatorname{sk}}}
\newcommand{\const}{{\operatorname{const}}}
\newcommand{\Sing}{{\operatorname{Sing}}}
\newcommand{\conj}{{\operatorname{Conj}}}
\newcommand{\Pic}{{\operatorname{Pic}}}
\newcommand{\Crit}{{\operatorname{Crit}}}
\newcommand{\Ch}{{\operatorname{Ch}}}
\newcommand{\discr}{{\operatorname{discr}}}
\newcommand{\Tor}{{\operatorname{Tor}}}
\newcommand{\Conj}{{\operatorname{Conj}}}
\newcommand{\val}{{\operatorname{val}}}
\newcommand{\Val}{{\operatorname{Val}}}
\newcommand{\res}{{\operatorname{res}}}
\newcommand{\add}{{\operatorname{add}}}
\newcommand{\tmu}{{\C\mu}}
\newcommand{\ov}{{\overline v}}\newcommand{\on}{{\overline n}}
\newcommand{\ox}{{\overline{x}}}
\newcommand{\tet}{{\theta}}
\newcommand{\Del}{{\Delta}}
\newcommand{\bet}{{\beta}}
\newcommand{\kap}{{\kappa}}
\newcommand{\del}{{\delta}}
\newcommand{\sig}{{\sigma}}
\newcommand{\alp}{{\alpha}}
\newcommand{\Sig}{{\Sigma}}
\newcommand{\Gam}{{\Gamma}}
\newcommand{\gam}{{\gamma}}
\newcommand{\Lam}{{\Lambda}}
\newcommand{\lam}{{\lambda}}
\newcommand{\SC}{{SC}}
\newcommand{\MC}{{MC}}
\newcommand{\nek}{{,...,}}
\newcommand{\cim}{{c_{\mbox{\rm im}}}}
\newcommand{\mathto}{\mathop{\to}}
\newcommand{\op}{{\overline p}}

\newcommand{\w}{{\omega}}

\title{New cases of logarithmic equivalence of Welschinger and
Gromov-Witten invariants}
\author{Ilia Itenberg \and Viatcheslav Kharlamov
\and Eugenii Shustin}
\date{}
\maketitle

{\hskip1in \it Dedicated to V.~I.~Arnol'd, at his 70th birthday
anniversary.}

\begin{figure}[htb]
\hfill\includegraphics[height=4cm,width=5cm,angle=0,draft=false]{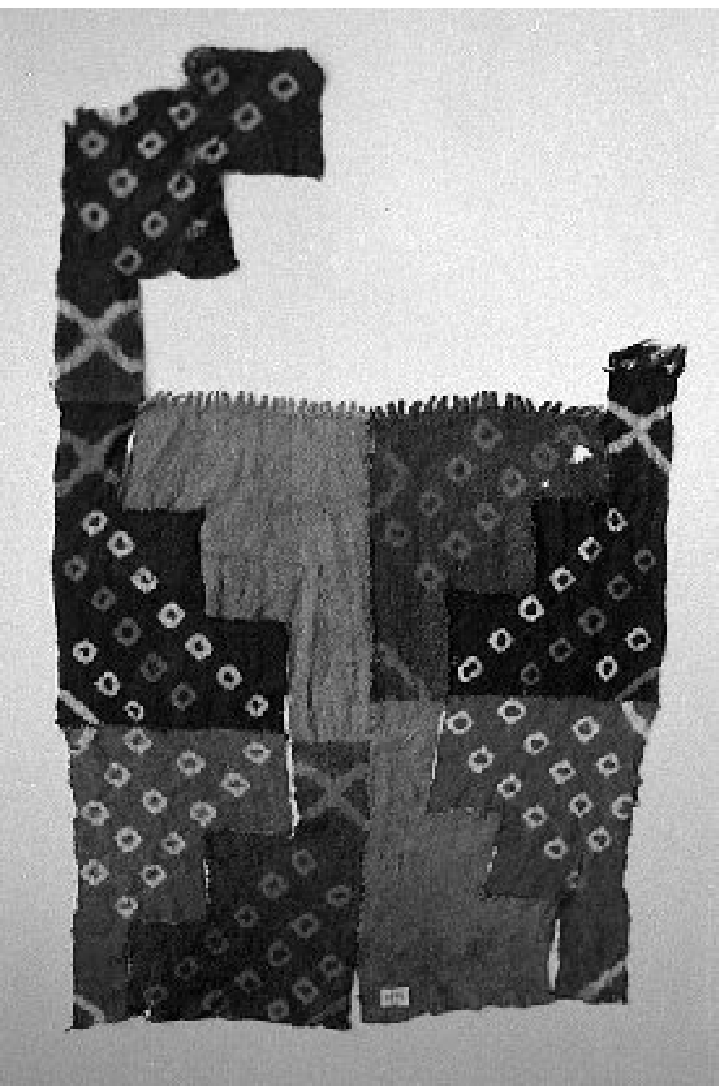}

{\hskip4.2in\it Nazca Tie-Dye Cloth}
\end{figure}

\begin{abstract} We consider
$\PP^1\times \PP^1$ equipped
with the complex conjugation $(x,y)\mapsto (\bar{y},\bar{x})$
and blown up in at most two real, or two complex conjugate, points.
For these four
surfaces we prove the logarithmic equivalence of Welschinger and
Gromov-Witten invariants.
\end{abstract}

\section{Introduction}\label{intro}
\footnotetext[1]{2000 {\it Mathematics Subject Classification}:
Primary 14N10. Secondary 14P99, 14N35.}
\footnotetext[2]{{\it Key words and phrases}:
Welschinger invariants, Gromov-Witten invariants, toric surfaces,
enumerative geometry, tropical curves.}
\footnotetext[3]{The authors were partially supported
by a grant from the Ministry of Science and Technology, Israel,
and
Minist\`{e}re des Affaires Etrang\`{e}res, France. The first two authors
were partially funded by the ANR-05-0053-01 grant of Agence
Nationale de la Recherche and a grant of Universit\'{e} Louis Pasteur,
Strasbourg.
The first two authors are participants of the CNRS-RFBR grant
"Problems in mathematical physics, tropical and idempotent
mathematics".
The third author acknowledges a support from the grant no. 465/04
from the Israel Science Foundation
and a support from the Hermann-Minkowski-Minerva Center
for Geometry at the Tel Aviv University.}

Welschinger invariants \cite{W, W1} applied to unnodal Del Pezzo
surfaces bound from below the number of real rational curves in a
given linear system
which pass
through a real generic collection of points.
In our previous papers~\cite{IKS1, IKS2}, using the methods of
tropical enumerative geometry developed by G.~Mikhalkin
\cite{Mi0,Mi} and E.~Shustin \cite{Sh0,Sh1}, we studied the toric
unnodal Del Pezzo surfaces with tautological real structure and
showed that for these surfaces Welschinger and Gromov-Witten
invariants are equivalent in the logarithmic scale, if all or
almost all fixed points in the generic collection are real. Here
we continue such an asymptotic study of Welschinger invariants and
consider non-tautological real structures on toric unnodal Del
Pezzo surfaces. Up to isomorphisms
respecting the real
structure, there are only five toric unnodal Del Pezzo surfaces
with a non-tautological real structure
and non-empty real part. One is obtained from
$\PP^1\times \PP^1$ equipped with the standard (tautological)
complex conjugation by blowing up two complex conjugate points,
and the four others are obtained from $\PP^1\times \PP^1$ equipped
with the complex conjugation $(x,y)\mapsto (\bar{y},\bar{x})$ by
blowing up at most two real, or two complex conjugate, points.

We
look at
collections of real points on
any of the four latter surfaces,
apply the tropical
formula elaborated in \cite{Sh2} to the multiples $nD$ of a real
ample divisor $D$ on
such a surface~$\Sigma$,
and
prove that Welschinger and Gromov-Witten invariants, $W_{\Sigma,nD}$ and
$GW_{\Sigma,nD}$,
are equivalent in the logarithmic scale: \mbox{$\log W_{\Sigma,nD}
=\log GW_{\Sigma,nD} +O(n)$}.
Recall that, as is shown in \cite{IKS2, IKS3}, \mbox{$\log
GW_{\Sigma,nD}
=(c_1(\Sigma)\cdot D) \, n\log {n} +O(n)$}.



\medskip

{\bf Acknowledgements}. A considerable part of this work was done
during our visits to the Max-Planck-Institut f\"{u}r Mathematik,
Bonn. We thank the MPIM for
hospitality and excellent working
conditions.

\section{Combinatorial bound}\label{invariants}

\begin{figure}
\begin{center}
\epsfxsize
95mm \epsfbox{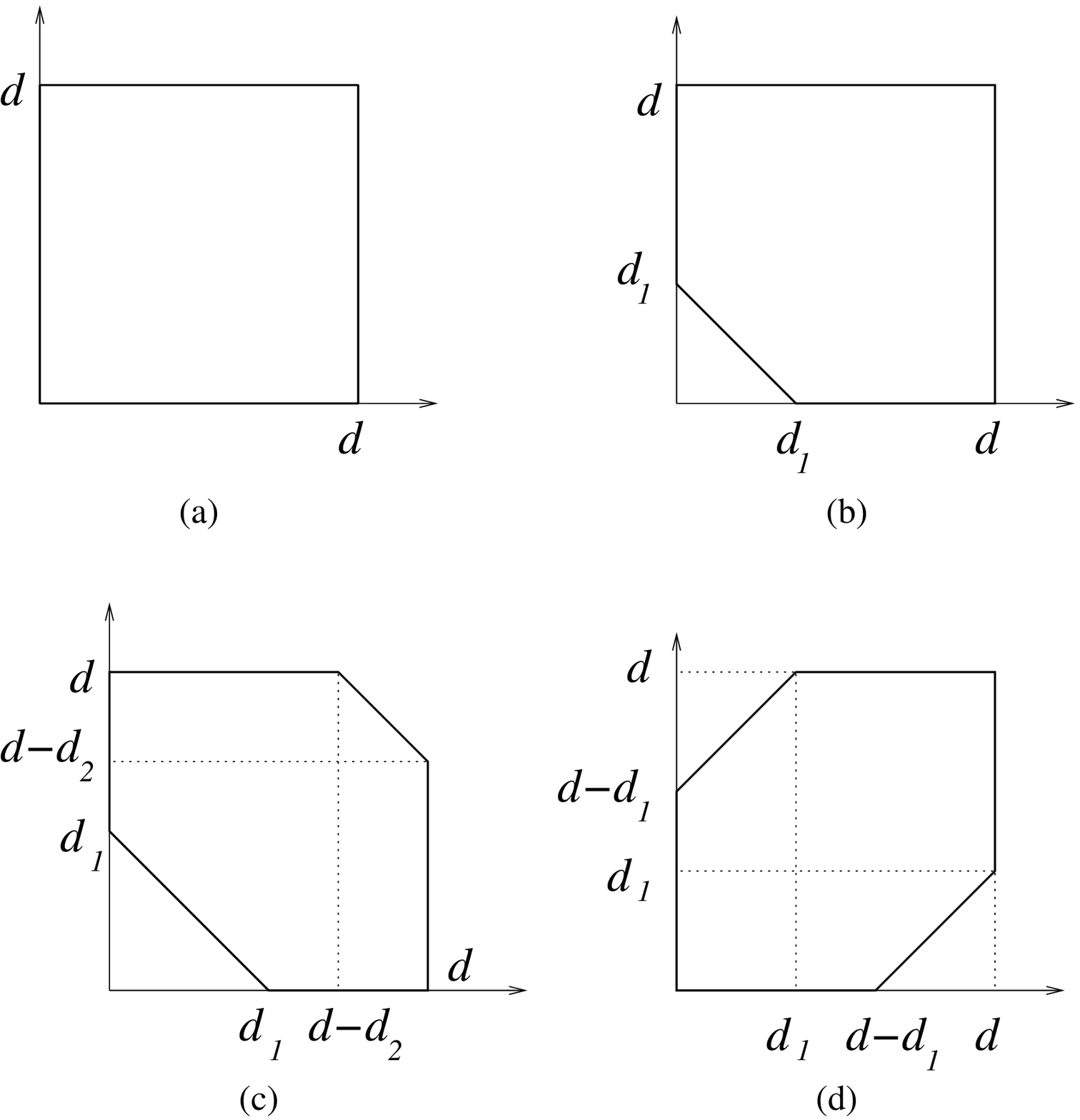}
\end{center}
\caption{Polygons defining $\bbS^2$, $\bbS^2_{1,0}$,
$\bbS^2_{2,0}$, and $\bbS^2_{0,2}$}
\label{fn2}
\end{figure}

As toric surfaces,
the four
real
Del Pezzo surfaces,
$\bbS^2$, $\bbS^2_{1,0}$, $\bbS^2_{2,0}$, and $\bbS^2_{0,2}$,
we deal with
are associated with the following convex
lattice polygons in $\R^2$ (see Figure \ref{fn2}):
\begin{itemize}
\item squares with vertices $(0,0),(d,0),(0,d),(d,d)$, where $d\ge 1$;
\item pentagons with vertices \hskip10pt
$(d,d),(0,d),(0,d_1),(d_1,0),(d,0)$, 
where $1\le d_1
<d$;
\item hexagons with vertices
$(0,d_1),(d_1,0),(d,0),(d,d-d_2), (d-d_2,d), (0,d)$,
where $1 \le d_2 \le d_1 < d$;
\item and hexagons with vertices
$(0,0),(d-d_1,0),(d,d_1),(d,d),(d_1,d),(0,d-d_1)$,
where $1\le d_1<d$.
\end{itemize}
For the toric surface~$\Sigma$ associated with such a polygon~$\Del$,
the real structure we study is the involution
which acts in
the principal orbit $(\C^*)^2 \subset \Sigma$
by $\conj(x,y)=(\overline y,\overline x)$.
Its natural lift to the
ample line bundle ${\cal L}_\Del$ 
generated by monomials $x^iy^j$, $(i,j)\in\Del\cap\Z^2$,
acts by
$\conj_*(a_{ij}x^iy^j)=\overline{a}_{ij}x^jy^i$, $(i,j)\in\Del$, 
and thus gives rise to
the reflection of~$\Del$ with respect to the bisectrix
${\cal B}$ of the positive quadrant.
Denote by~$D_\Del$ (or simply $D$) an ample divisor
which defines ${\cal L}_\Del$.

The goal of this section is to deduce from~\cite{Sh2}, Theorem~1.1,
a lower bound for the Welschinger invariant~$W_{\Sigma,D}$.
\footnote{In~\cite{Sh2}, this invariant is denoted by $W_0(\Sigma, {\cal L})$,
where ${\cal L} = {\cal L}_\Del$.}
To this end, we introduce the following objects.

For each integer point belonging to the boundary of~$\Del$,
trace the straight
line
through
this point and its
image under the reflection with respect to~$\cal B$.
The union of all the traced lines cuts
${\cal B} \cap \Del$ in certain segments; denote their number
by~$m$. Identify ${\cal B} \cap \Del$ with the segment $[0, m]
\subset \R$ in such a way that the intersection points of ${\cal
B} \cap \Del$ with the traced lines are mapped to the integer
points of $[0,m]$. To each integer point $i \in [0,m]$ associate a
non-negative integer number $\sig(i)$ equal to the
integer
length of the intersection of
the corresponding straight line with~$\Del$.

A finite multi-set of closed intervals in $\R$
is called a $\Del$-{\it proper
system} (or simply {\it proper system}) if
\begin{itemize}
\item each interval is contained in $[0,m]$ and
has integer endpoints (intervals
reduced to a
point are allowed),
\item the total number of intervals is $|\partial\Del| - m - 1$,
where $|\partial\Del|$ is the integer length of the boundary of~$\Del$,
\item for any integer $i \in [0,m]$, the number of intervals
containing~$i$ is equal to $\sig(i)$.
\end{itemize}

Given a $\Del$-proper system, consider
the disjoint union~$g'$ of the intervals of the system,
and complete~$g'$ to a graph~$g$ introducing
$m$ additional vertices indexed by the half-integer points
$i + 1/2$, $i = 0, \ldots , m - 1$,
and additional edges connecting each point $i + 1/2$
with all the right endpoints $i$ and all the left endpoints $i + 1$
of the intervals in~$g'$.

A $\Del$-proper system is called {\it admissible},
if its graph~$g$ is a tree. An admissible
$\Del$-proper system is {\it marked},
if it is equipped with a {\it marking} which associates to each
interval~$I$ of the system an integer point of~$I$;
the latter point is called {\it marked}.

The following statement
is an immediate consequence of~\cite{Sh2}, Theorem~1.1.

\begin{lemma}\label{the-lemma}
Let~$\Del$ be one of the polygons shown in Figure~\ref{fn2},
and $\Sig$ the toric surface
associated with~$\Del$ and equipped
with the real structure~$\conj$ {\rm (}described above{\rm )}.
Then, the Welschinger invariant $W_{\Sig,D_\Del}$ is greater
than or equal to the number of marked admissible $\Del$-proper
systems.
\end{lemma}

\section{Logarithmic asymptotics}\label{asymptotics}

\subsection{Main theorem}\label{main-theorem}

\begin{theorem}\label{tsp2}
Let~$\Sig$ be one of the real surfaces
$\bbS^2$, $\bbS^2_{1,0}$, $\bbS^2_{2,0}$, and
$\bbS^2_{0,2}$.
For any real ample divisor~$D$ on~$\Sig$,
it holds
\begin{equation}
\log W_{\Sig, nD} = (c_1(\Sig)\cdot D)n\log n + O(n).
\end{equation}
In particular,
\begin{equation}
\lim_{n\to\infty}\frac{\log W_{\Sig, nD}}{\log GW_{\Sig, nD}}= 1,
\label{spe200}\end{equation}
where $GW_{\Sig, nD}$ is the genus zero
Gromov-Witten
invariant.
\end{theorem}


Since $W_{\Sig, nD} \leq GW_{\Sig, nD}$ and
$GW_{\Sig, nD} = (c_1(\Sig)\cdot D)n\log n + O(n)$,
to prove Theorem~\ref{tsp2} it is sufficient
to prove the lower bound $W_{\Sig, nD} \geq (c_1(\Sig)\cdot D)n\log n + O(n)$.
Due to Lemma~\ref{the-lemma}
and the identity $|\partial\Del|=c_1(\Sig)\cdot D$, the latter lower bound
would follow from the inequality
\begin{equation}
\log S_{n\Del} \geq |\partial\Del|\cdot n \log
n+O(n),\label{spe201}\end{equation}
where $S_{n\Del}$ is the number of marked admissible $n\Del$-proper
systems.
This inequality is proved in Sections~\ref{0-points}, \ref{(1,0)},
\ref{(2,0)}, and~\ref{(0,2)},
where each of the surfaces
$\bbS^2$, $\bbS^2_{1,0}$, $\bbS^2_{2,0}$, and $\bbS^2_{0,2}$
is treated separately.


\subsection{Admissibility}\label{admissibility}

Let~$\Gam$ be a finite set of disjoint horizontal segments
with integer endpoints in $\R^2$ (degenerated segments are allowed).
For any vertical strip $b = \{i \leq x \leq i+1\}$,
where~$i$ is an integer, denote by
$\Gam^L(b)$ (respectively, $\Gam^R(b)$)
the subset of~$\Gam$ formed by the segments
whose right endpoint belongs to $x = i$
(respectively, left endpoint belongs to $x = i + 1$).

\begin{lemma}\label{sink}
Assume that~$\Gam$ can be represented
as the disjoint union of two subsets~$\Gam_L$ and~$\Gam_R$
satisfying the following properties:
\begin{enumerate}
\item[(i)] for any vertical strip $b = \{i \leq x \leq i+1\}$
such that~$i$ is an integer,
the union of $\Gam_R \cap \Gam^R(b)$
and $\Gam_L \cap \Gam^L(b)$ contains at most one element,
\item[(ii)] if the union of $\Gam_R \cap \Gam^R(b)$
and $\Gam_L \cap \Gam^L(b)$ contains an element~$s$,
no element of~$\Gam^L(b) \cup \Gam^R(b)$ lies below~$s$;
\item[(iii)] there exists exactly
one vertical strip $b = \{i \leq x \leq i+1\}$
such that~$i$ is an integer, at least one of the sets~$\Gam^L(b)$
and~$\Gam^R(b)$ is nonempty, and the union of $\Gam_R \cap \Gam^R(b)$
and $\Gam_L \cap \Gam^L(b)$ is empty.
\end{enumerate}
If the projections of segments of~$\Gam$ on the horizontal axis
form a proper system, then this proper system is admissible.
\end{lemma}

{\bf Proof.}
For a proper system as in the lemma,
identify~$\Gam$ with the disjoint union $g'$ of the intervals of the system,
and consider the graph~$g$ as in Section~\ref{invariants}.
Orient the segments of~$\Gam_L$ to the left,
the segments of~$\Gam_R$ to the right,
and orient each additional edge of~$g$
by extending the orientation of the adjacent horizontal segment.
The conditions~($i$) and~($ii$) give a deformation retraction
of~$g$ to a finite set of vertices, and the condition~($iii$)
guarantees that the latter set has only one element.
\proofend

\subsection{Case $\Sig=\bbS^2$}\label{0-points}

Let~$\Del$ be the square shown in Figure~\ref{fn2}(a).
In this case, the required inequality~(\ref{spe201})
reads as
\begin{equation}
\log S_{n\Del}\ge4dn\log n + O(n)\
.\label{spe203-1}
\end{equation}

\begin{figure}
\begin{center}
\epsfxsize 145mm \epsfbox{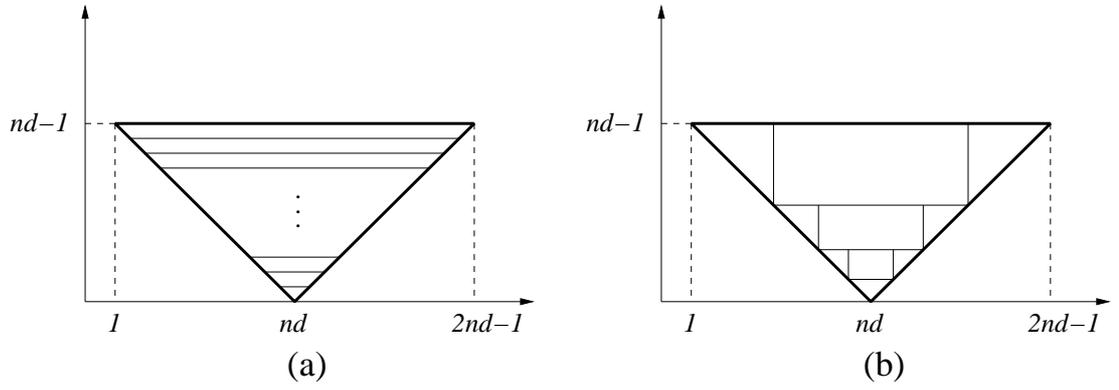}
\end{center}
\caption{First steps
in the construction of admissible systems for $\bbS^2$}
\label{triangle}
\end{figure}

To construct an appropriate number of marked admissible
$n\Del$-proper systems, consider the triangle~$T(n,d)$ with
vertices $(1, nd-1)$, $(nd, 0)$, and \mbox{$(2nd-1, nd-1)$} (see
Figure~\ref{triangle}(a)). At each integer level $y = j$, $0 \leq
j \leq nd-1$ consider the maximal horizontal segment contained
in~$T(n,d)$.
If $j\ne0$,
make a {\it hole} in the considered segment
by removing an open unit
interval
with integer endpoints.
This {\it perforation procedure} gives rise to a set of $2nd - 1$
horizontal segments whose projections form an $n\Del$-proper
system.

Inscribe in $T(n,d)$ a sequence of maximal size rectangles~$R_i$
satisfying the following properties: each rectangle is symmetric
with respect to the vertical line $x = nd$, and the length of
horizontal edges of each rectangle is twice the length of its
vertical edges (see Figure~\ref{triangle}(b)). The right upper
vertices $(x_i, y_i),\ i\ge 1$ of these rectangles are given by
$$
x_1= nd + \left[\frac{nd-1}{2}\right],\
y_1 = nd - 1,\
y_{i+1}=y_i - \left[\frac{y_i}{2}\right],\ x_{i+1} = nd+y_{i+1}-
\left[\frac{y_{i+1}}{2}\right].
$$
Let~$k$ be the number of rectangles. Notice that $y_k = 2$,
and put $y_{k+1}=y_k - \left[\frac{y_k}{2}\right]=1$.

Restrict
the choice of holes
in the perforation procedure in the following way:
\begin{itemize}
\item all the holes are contained in the half-plane $x \geq nd$,
\item for any integer $1 \leq i \leq k$ all the holes at
the levels $y_{i+1} + 1 \leq y \leq y_i$
are contained in~$R_i$,
\item for any integer $1 \leq i \leq k - 1$ no two holes
at the levels $y_{i+1} + 1 \leq y \leq y_i$ have the same projection
on the horizontal axis.
\end{itemize}
The set of segments obtained via such a
perforation procedure is called a {\it perforated} $(n,d)$-{\it collection}.
The number $M(n,d)$ of
perforated $(n,d)$-collections
is
equal to
$$(y_1 - y_2)!(y_2 - y_3)! \ldots (y_k - y_{k+1})! \ .$$
According to the Styrling formula,
$$\log M(n,d)
= ((y_1 - y_2) + (y_2 - y_3) + \ldots + (y_k - y_{k+1})) \log n + O(n)
= dn \log n + O(n).$$

\begin{figure}
\begin{center}
\epsfxsize 105mm \epsfbox{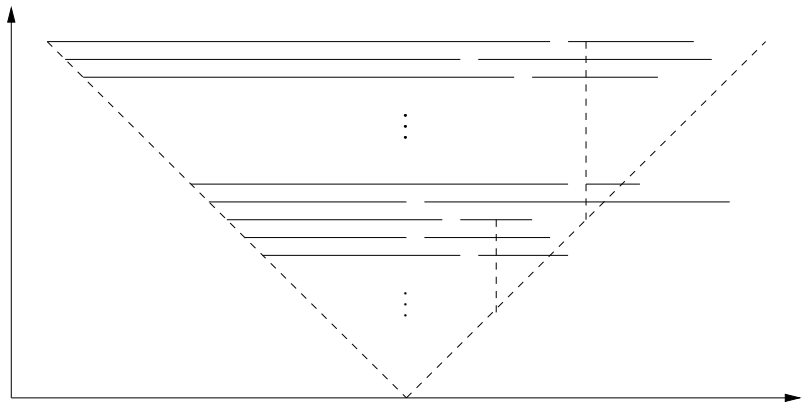}
\end{center}
\caption{A permuted perforated $(n,d)$-collection}
\label{permutation}
\end{figure}

For any
perforated $(n,d)$-collection~$c$ and any permutations $\sig_1$,
$\ldots$, $\sig_{k-1}$, where $\sig_i$, $i = 1$, $\ldots$, $k -
1$, is a permutation of $\{y_{i+1} + 1, y_{i+1} + 2, \ldots ,
y_i\}$, consider the set of segments $c_{\sig_1, \ldots , \sig_{k
- 1}}$ obtained from~$c$ in the following way: for each integer $1
\leq i \leq k -1$, cut along the vertical line $x = x_i$ the
segments of~$c$ lying on the levels $y_{i+1} + 1 \leq y \leq y_i$
and intersecting the line $x = x_i$, permute according to~$\sig_i$
the right-hand parts of
the segments we have cut,
and glue the adjacent parts
in order to form new segments
(see Figure~\ref{permutation}).
The set
$c_{\sig_1, \ldots , \sig_{k - 1}}$ is called a {\it permuted
perforated} $(n,d)$-{\it collection}. It consists of the point
$(nd, 0)$ and two segments at each integer level $1 \leq y \leq nd
- 1$.

The number ${\widetilde M}(n,d)$
of the permuted perforated $(n,d)$-collections
$c_{\sig_1, \ldots , \sig_{k-1}}$,
where~$c$ runs over all the perforated $(n,d)$-collections
and $\sig_i$, $i = 1$, $\ldots$, $k - 1$, runs over
all the permutations of $\{y_{i + 1} + 1, y_{i + 1} + 2, \ldots y_i\}$,
is equal to
$$M(n,d)(y_1 - y_2)!(y_2 - y_3)! \ldots (y_k - y_{k+1})! \ .$$
Thus,
$\log {\widetilde M}(n,d) = 2dn \log n + O(n)$.

The projection on the horizontal axis
of any permuted perforated $(n,d)$-collection
is an $n\Del$-proper system.
The restriction imposed above on the choice of holes
guarantees that the projection
of all the permuted perforated $(n,d)$-collections
produces ${\widetilde M}(n,d)$
pairwise distinct $n\Del$-proper systems.
All the resulting systems
are admissible as it follows from Lemma~\ref{sink}
applied to any permuted perforated $(n,d)$-collection represented
as the disjoint union of the segments lying on the left-hand side
of the holes (the subset~$\Gamma_R$) and the segments lying on the
right-hand side of the holes (the subset $\Gamma_L$).

Mark each
of ${\widetilde M}(n,d)$
admissible $n\Del$-proper systems as above
in such a way that
\begin{itemize}
\item no marked point of the projection of a segment at level~$1$
does coincide with the point~$nd$,
\item for any integer~$i$ between~$1$ and~$k$,
the marked points of
the projections of segments at any level $y_{i+1} + 1 \leq y \leq y_i$
are placed outside of
the projection of~$R_i$.
\end{itemize}

For each system, this can be done in
$\left((y_1 - y_2)!(y_2 - y_3)! \ldots (y_k - y_{k+1})!\right)^2$ different
ways. Thus, the logarithm of the number of obtained marked
admissible $n\Del$-proper systems is $4dn \log n + O(n)$.
This proves Theorem~\ref{tsp2} in the case $\Sig = \bbS^2$.


\subsection{Case $\Sig=\bbS^2_{1,0}$}\label{(1,0)}

Let~$\Del$ be the pentagon shown in Figure~\ref{fn2}(b).
In this case, the required inequality~(\ref{spe201})
reads as
\begin{equation}\log S_{n\Del}\ge(4d-d_1)n\log n + O(n)\
.\label{spe203-2}
\end{equation}

\begin{figure}
\begin{center}
\epsfxsize 145mm \epsfbox{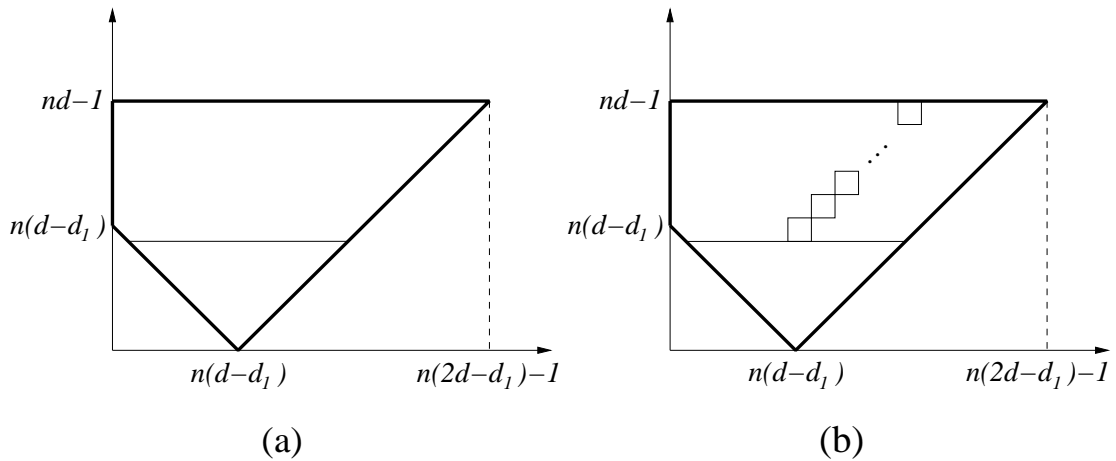}
\end{center}
\caption{Construction of admissible systems for $\bbS^2_{1,0}$}
\label{quadrangle}
\end{figure}

To construct an appropriate number of marked admissible
$n\Del$-proper systems, consider the quadrangle~$Q(n,d,d_1)$ with
vertices
$$(0, nd-1), \ (0, n(d-d_1)), \ (n(d-d_1), 0), \
\text{\rm and} \ (n(2d - d_1)-1, nd-1),$$ (see
Figure~\ref{quadrangle}(a)). For the triangle $T(n,d-d_1) \subset
Q(n,d,d_1)$ use the construction described in
Section~\ref{0-points}. To complete
the resulting permuted perforated $(n,d-d_1)$-collections,
we proceed in the following way.

Consider the up-right staircase~$E$ formed by squares of size $n
\times n$ such that~$E$ starts at the middle point $(n(d - d_1),
n(d - d_1) - 1)$ of the upper side of~$T(n, d - d_1)$ (see
Figure~\ref{quadrangle}(b)). At each integer level $y = j$, $n(d -
d_1 )\leq j \leq nd-1$ consider the maximal horizontal segment
contained in~$Q(n, d, d_1)$, and use the perforation procedure
(that is, make a hole in each segment considered) choosing holes
in such a way that all these holes are contained in~$E$, no hole
is taken on the lower sides of the squares forming~$E$, and no two
holes have the same projection on the horizontal axis. This gives
$(n!)^{d_1}$ sets of segments. For any of these sets and any
permuted perforated $(n,d-d_1)$-collection, their union is called
a {\it perforated} $(n,d,d_1)$-{\it collection}.

The projection to the horizontal axis of any
perforated $(n,d,d_1)$-collection is an $n\Del$-proper system.
Due to Lemma~\ref{sink}, any resulting $n\Del$-proper system
is admissible. For any such system, there are at least
$(nd_1)!(nd_1)!$ choices of marking for the projections
of segments lying above $T(n, d-d_1)$.
Thus, the logarithm of the number of marked admissible
$n\Del$-proper systems is at least
$$4(d - d_1)n \log n + O(n) + d_1\log n! + 2 \log (nd_1)!
= (4d - d_1)n \log n + O(n).$$ This proves Theorem~\ref{tsp2} in
the case $\Sig = \bbS^2_{1,0}$.


\subsection{Case $\Sig=\bbS^2_{2,0}$}\label{(2,0)}

Let~$\Del$ be the hexagon shown in Figure~\ref{fn2}(c).
In this case, the required inequality~(\ref{spe201})
reads as
\begin{equation}\log S_{n\Del}\ge(4d-d_1-d_2)n\log n + O(n)\
.\label{spe203-3}
\end{equation}

\begin{figure}
\begin{center}
\epsfxsize 145mm \epsfbox{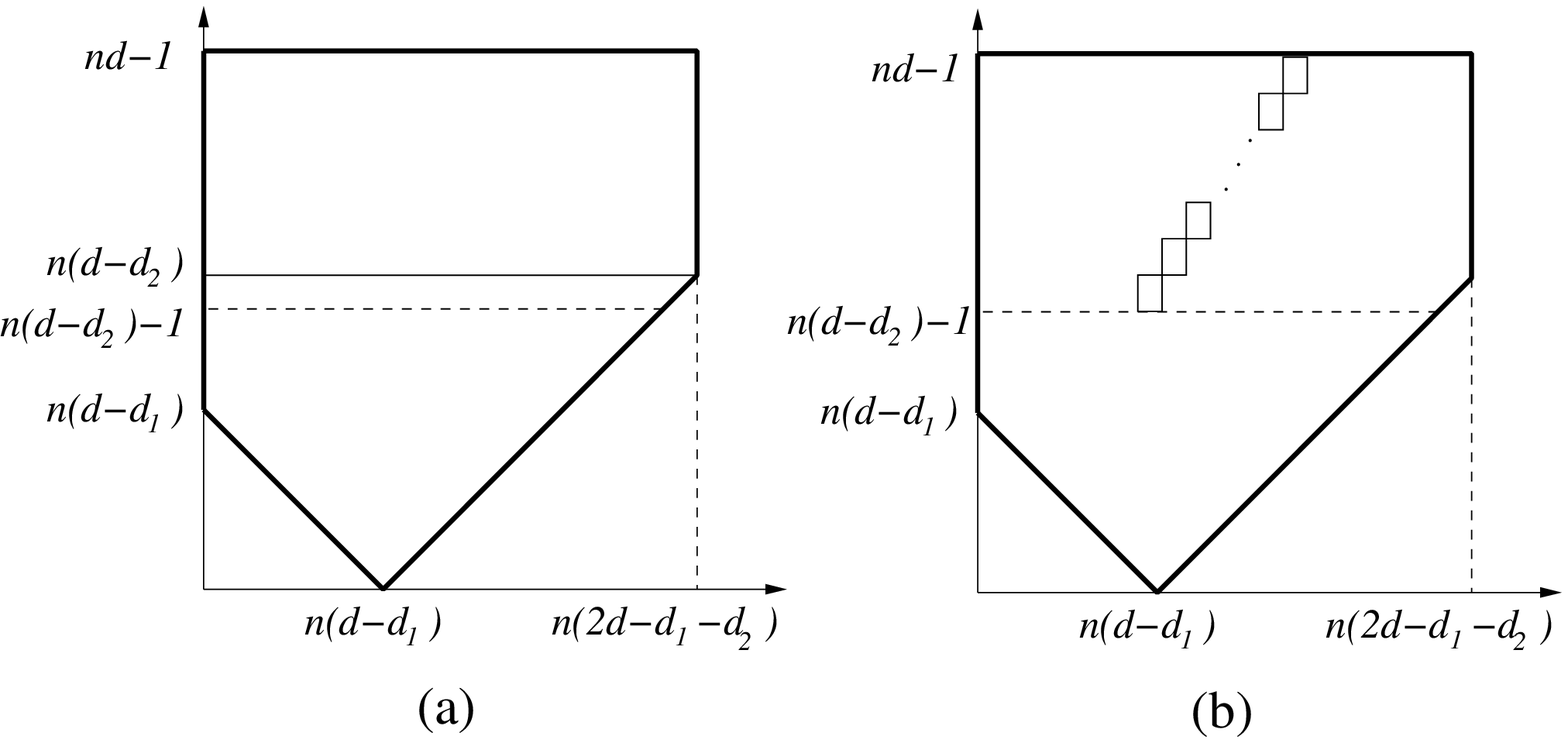}
\end{center}
\caption{Construction of admissible systems for $\bbS^2_{2,0}$}
\label{pentagon}
\end{figure}

To construct an appropriate number of marked admissible
$n\Del$-proper systems, consider the pentagon~$P(n,d,d_1,d_2)$ with
vertices
$$\displaylines{
(0, nd-1), \ (0, n(d-d_1)), \ (n(d-d_1), 0), \cr (n(2d - d_1 -
d_2), n(d - d_2))), \ \text{\rm and} \ (n(2d - d_1 - d_2), nd-1),}$$
(see Figure~\ref{pentagon}(a)).
For the quadrangle
$Q(n,d-d_2,d-d_1) \subset P(n,d,d_1,d_2)$ use the construction
described in Section~\ref{(1,0)}. To complete
the resulting perforated $(n,d-d_2,d_1-d_2)$-collections,
we proceed in the following way.

The remaining part of $P(n, d, d_1, d_2)$
is formed by a horizontal strip of height~$1$
and a rectangle of width~$n(2d - d_1 - d_2)$
and height~$nd_2 - 1$, see~Figure~\ref{pentagon}(a).
Consider an up-right staircase
\begin{itemize}
\item starting
at a point $(x_0, n(d - d_2) - 1)$
with $x_0 \geq [n(2d-d_1-d_2)/4]$,
\item ending
at a point $(x_1, nd - 1))$ with $x_1 \leq [3n(2d-d_1-d_2)/4] - 1$,
\item and formed by rectangles such that each rectangle
is of width~$1$ and of positive height smaller than or equal to $a
= \left[\frac{2d_2}{2d-d_1-d_2}\right] + 1$,
\end{itemize}
(see Figure~\ref{pentagon}(b)).
At each integer level $y = j$,
$n(d - d_2) \leq j \leq nd-1$ consider the maximal horizontal
segment contained in~$P(n, d, d_1, d_2)$, and use the perforation
procedure choosing holes in the rectangles of the staircase in
such a way that no hole is taken on the lower sides of the
rectangles. For any perforated $(n,d-d_2,d_1 - d_2)$-collection,
its union with the constructed set of segments is called a {\it
perforated} $(n,d,d_1,d_2)$-{\it collection}.

The projection to the horizontal axis of any perforated
$(n,d,d_1,d_2)$-collection is an $n\Del$-proper system. Due to
Lemma~\ref{sink}, any resulting $n\Del$-proper system is
admissible. For any such system, there are at least $[n(2d - d_1 -
d_2)/4]^{2nd_2}(a!)^{-2nd_2}$ choices of marking for the
projections of segments lying above $Q(n, d-d_2, d_1 - d_2)$.
Thus, the logarithm of the number of marked admissible
$n\Del$-proper systems is at least
$$(4(d - d_2) - (d_1 - d_2))n \log n + O(n) + 2d_2n\log n + O(n)
= (4d - d_1 - d_2)n \log n + O(n).$$ This proves
Theorem~\ref{tsp2} in the case $\Sig = \bbS^2_{2,0}$.

%

\subsection{Case $\Sig=\bbS^2_{0,2}$}\label{(0,2)}

Let~$\Del$ be the hexagon shown in Figure~\ref{fn2}(d).
In this case, the required inequality~(\ref{spe201})
reads as
\begin{equation}\log S_{n\Del}\ge(4d-2d_1)n\log n + O(n)\
.\label{spe203-4}
\end{equation}

\begin{figure}
\begin{center}
\epsfxsize 105mm \epsfbox{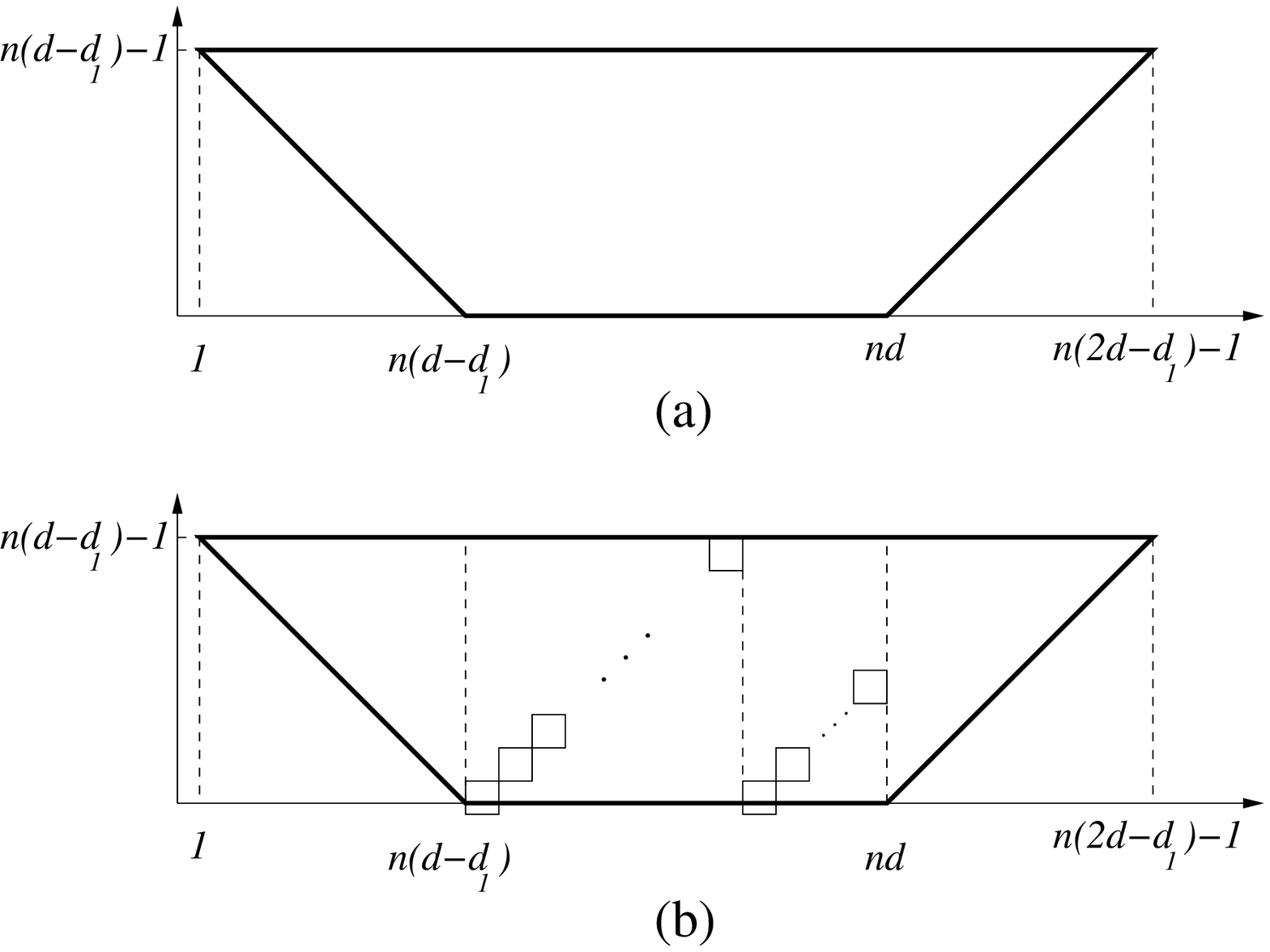}
\end{center}
\caption{Construction of admissible systems for $\bbS^2_{0,2}$}
\label{trapeze}
\end{figure}

To construct an appropriate number of marked admissible
$n\Del$-proper systems, consider the trapeze~$K(n,d,d_1)$ with
vertices
$$(1, n(d-d_1)-1, \ (n(d-d_1),0), \ (nd, 0), \
\text{\rm and} \ (n(2d - d_1)-1, n(d-d_1)-1),$$ (see
Figure~\ref{trapeze}(a)).

Consider the sequence~$C$ of up-right staircases
formed by squares of size $n \times n$
such that
\begin{itemize}
\item all the staircases of~$C$ are contained
in the vertical strip
$${\cal B} = \{n(d - d_1) \leq x \leq nd\},$$
\item each staircase starts at the level $y = -1$,
\item each staircase ends at the upper side of $K(n, d, d_1)$,
the only possible
exception being the last staircase,
\item the first staircase starts at
the point $(-1, n(d - d_1))$,
\item for each staircase, except the first one, the vertical line
where the staircase starts coincides with the vertical line
where the preceding staircase ends,
\end{itemize}
(see Figure~\ref{trapeze}(b); in the case $d_1 \leq d - d_1$, there
is only one staircase in~$C$). At each integer level $y = j$, $0
\leq j \leq n(d - d_1) - 1$ consider the maximal horizontal
segment contained in~${\cal B}$,
and use the perforation procedure
(this time we authorize several holes at the same level)
by choosing holes in such a way that all these holes are contained
in~$C$,
no hole is taken on the lower sides of the squares
forming the staircases,
and there is exactly one hole in each integer vertical strip
$i \leq x \leq i + 1$
contained in~${\cal B}$.
This gives $(n!)^{d_1}$ sets of segments.

Pick a permuted perforated $(n,d-d_1)$-collection~$\pi$
in $T(n,d-d_1)$, cut~$\pi$ along the vertical line $x = n(d - d_1)$,
keep the left half of~$\pi$ at its place
and shift the right half by the vector $(nd_1, 0)$.
The result of gluing of the obtained collection
with a set of segments constructed in ${\cal B}$
as described above
is called a perforated $K(n, d, d_1)$-collection.
The projection to the horizontal axis
of any perforated $K(n, d, d_1)$-collection
is a $n\Del$-proper system.

Any resulting $n\Del$-proper system is admissible.
Indeed, let~$\gamma$ be a perforated $K(n, d, d_1)$-collection.
Identifying~$\gamma$ with the disjoint union $g'$ of the intervals
of the projection of~$\gamma$ to the horizontal axis,
consider the graph~$g$ as in Section~\ref{invariants}.
In each integer vertical strip
$i \leq x \leq i + 1$
contained in~${\cal B}$, there is exactly one pair
of additional edges of~$g$, and this pair
fill up the only hole in $i \leq x \leq i + 1$.
Once the holes in~${\cal B}$ are filled up, Lemma~\ref{sink}
applies. This proves the admissibility of the projection of~$\gamma$.

Consider a perforated $K(n, d, d_1)$-collection~$\gamma$ obtained
by gluing of a permuted perforated $(n, d-d_1)$-collection~$\pi$
with a set of segments constructed in~${\cal B}$ as above. Any
marking of the projection of~$\pi$ can be extended to a marking of
the projection of~$\gamma$ via a choice of an integer point on
each segment entering under a staircase. The latter choice can be
done in at least $(nb_1)! \ldots (nb_k)!$ ways,
where $b_1$, $\ldots$, $b_k$ are the numbers of stairs in the
staircases (in fact, $b_1 =\ldots = b_{k-1}$).
Thus, the logarithm of the number of marked admissible
$n\Del$-proper systems is at least
$$\displaylines{
4(d - d_1)n \log n + O(n) + d_1\log n! + n(b_1 + \ldots + b_k)\log
n + O(n) \cr = (4d - 2d_1)n \log n + O(n).}$$ This proves
Theorem~\ref{tsp2} in the case $\Sig = \bbS^2_{0,2}$.

%

{\ncsc Universit\'{e} Louis Pasteur et IRMA \\[-21pt]

7, rue Ren\'{e} Descartes, 67084 Strasbourg Cedex, France} \\[-21pt]

{\it E-mail address}: {\ntt itenberg@math.u-strasbg.fr}

\vskip10pt

{\ncsc Universit\'{e} Louis Pasteur et IRMA \\[-21pt]

7, rue Ren\'{e} Descartes, 67084 Strasbourg Cedex, France} \\[-21pt]

{\it E-mail address}: {\ntt kharlam@math.u-strasbg.fr}

\vskip10pt

{\ncsc School of Mathematical Sciences \\[-21pt]

Raymond and Beverly Sackler Faculty of Exact Sciences\\[-21pt]

Tel Aviv University,
Ramat Aviv, 69978 Tel Aviv, Israel} \\[-21pt]

{\it E-mail address}: {\ntt shustin@post.tau.ac.il}

\end{document}